\documentclass[10pt,twoside,reqno]{amsart}
\usepackage{amssymb}
\calclayout

\numberwithin{equation}{section}
\makeatletter

\renewcommand{\@secnumfont}{\bfseries}

\renewcommand{\section}{\@startsection{section}{1}%
  {0mm}{.7\linespacing\@plus\linespacing}{.5\linespacing}
  {\normalfont\bfseries\centering}}

\newcommand{\bibsection}{\@startsection{section}{1}%
  {0mm}{.7\linespacing\@plus\linespacing}{.5\linespacing}
  {\normalfont\scshape\centering}}

\renewcommand{\@biblabel}[1]{#1.}

\newtheorem{thm}{\bf Theorem}[section]

\newtheorem{cor}[thm]{\bf Corollary}

\begin{document}

\vspace{1.3cm}

\title  {A note on $q$-analogue of Boole polynomials}

\author{ Dae San Kim, TAEKYUN KIM, JONG JIN SEO}
\thanks{\scriptsize }

\address{1\\Department of  Mathematics\\
            Sogang University\\
           Seoul 121-742, Republic of Korea}
\email{dskim@kw.ac.kr}

\address{2\\Department of  Mathematics\\
            Kwangwoon National University\\
           Seoul 139-701, Republic of Korea}
\email{tkkim@kw.ac.kr}

\address{3\\Department of Applied Mathematics\\
            Pukyong National University\\
            Busan 608-737, Republic of Korea}
\email{seo2011@pknu.ac.kr}


\maketitle

\begin{abstract} In this paper, we consider the $q$-extensions of Boole polynomials. From those polynomials, we derive some new and interesting properties and identities related to special polynomials.
\end{abstract}

\pagestyle{myheadings}
\markboth{\centerline{\scriptsize  D. S. Kim, T. Kim,  J.J. Seo}}
          {\centerline{\scriptsize A note on $q$-analogue of Boole polynomials}}
\bigskip
\bigskip
\medskip
\section{\bf Introduction}
\bigskip
\medskip

Let $p$ be a fixed odd prime number. Throughout this paper, $ \mathbb{Z}_p, \mathbb{Q}_p$ and $\mathbb{C}_p$ will denote the ring of $p$-adic integers, the field of $p$-adic numbers and the completion of algebraic closure of $\mathbb{Q}_p$. The $p$-adic norm ${|\cdot|}_{p}$ is normalized as $|p|_{p}={1/p}$. The space of continuous functions on $\mathbb{Z}_p$ is denoted by $C\left(\mathbb{Z}_p\right)$. Let $q$ be an indeterminate in $\mathbb{C}_p$ with $|1-q|_p<p^{-1/p-1}$. The $q$-number of $x$ is defined by $[x]_{q}={{1-q^x}\over{1-q}}$. Note that $\lim_{q\rightarrow 1} {[x]_q=x}$. For $f \in C\left(\mathbb{Z}_p\right)$, the fermionic $p$-adic $q$-integral on $\mathbb{Z}_p$ is defined by Kim to be
\begin{equation}\begin{split}
&I_{-q}(f)=\int_{\mathbb{Z}_p} f(x) d\mu_{-q} (x) = \lim_{N\rightarrow \infty} {1\over {[p^N]_{-q}}} \sum_{x=0}^{p^N-1} f(x)(-1)^{x},\\
& \textnormal{where} \ [x]_{-q}={{1-(-q)^{x}}\over {1+q}}\ (\textnormal{see}\ [1-9]).
\end{split}\end{equation}
From (1.1), we note that
\begin{equation}\begin{split}
&q^{n}I_{-q}(f_{n})+(-1)^{n-1}I_{-q}(f)=[2]_{q}\sum_{l=0}^{n-1}(-1)^{n-1-l}q^{l}f(l),\\
&\textnormal{where} \ f_{n}(x)=f(x+n), (n\geq1)\ (\textnormal{see}\ [4]).
\end{split}\end{equation}
In particular, for $n$=1,
\begin{equation}
qI_{-q}(f_{1})+I_{-q}(f)=[2]_{q}f(0).
\end{equation}
As is well known, the Boole polynomials are defined by the generating function to be
\begin{equation}
\sum_{n=0}^{\infty}Bl_{n}(x|\lambda){t^{n}\over{n!}}={1\over{1+(1+t)^{\lambda}}}(1+t)^{x}, \ (\textnormal{see}\ [2,11]).
\end{equation}
When $\lambda=1, 2Bl_{n}(x|1)=Ch_{n}(x)$ are Changhee polynomials which are defined by
\begin{equation}
{2\over{t+2}}(1+t)^{x}=\sum_{n=0}^{\infty}Ch_{n}(x){t^{n}\over{n!}}\ (\textnormal{see}\ [2]).
\end{equation}
The Euler polynomials of order $\alpha$ are defined by the generating function to be
\begin{equation}
\left({2\over{e^{t}+1}}\right)^{\alpha}e^{xt}=\sum_{n=0}^{\infty} E_{n}^{(\alpha)}(x){t^{n}\over{n!}}, \ (\textnormal{see}\ [2,11]).
\end{equation}
When $x=0, E_{n}^{(\alpha)}=E_{n}^{(\alpha)}(0)$ are called the Euler numbers of order $\alpha$.\\
In particular, for $\alpha=1, E_{n}(x)=E_{n}^{(1)}(x)$ are called the ordinary Euler polynomials.\\
The Stirling number of the first kind is given by the generating function to be
\begin{equation}
\log{(1+t)^{m}}=m!\sum_{l=m}^{\infty}S_{1}(l,m){t^{l}\over{l!}},(m\geq0),
\end{equation}
and the Stirling number of the second kind is defined by the generating function to be
\begin{equation}
(e^{t}-1)^{m}=m!\sum_{l=m}^{\infty}S_{2}(l,m){t^{l}\over{l!}},\ (\textnormal{see}\ [11]).
\end{equation}
In this paper, we consider the $q$-extensions of Boole polynomials. From those polynomials, we derive new and interesting properties and identities related to special polynomials.

\section{\bf $q$-analogue of Boole polynomials}
\bigskip
\medskip
In this section, we assume that $t\in\mathbb{C}_p$ with ${|t|}_{p}<p^{-1\over{p-1}}$ and $\lambda\in\mathbb{Z}_p$ with $\lambda\neq0$. From (1.3), we note that
\begin{equation}\begin{split}
\int_{\mathbb{Z}_p}(1+t)^{x+\lambda y}d\mu_{-q}(y) &={1+q\over{1+q(1+t)^{\lambda}}}(1+t)^{x}\\
&=\sum_{n=0}^{\infty}[2]_{q}Bl_{n,q}(x|\lambda){t^{n}\over{n!}},
\end{split}\end{equation}
where $Bl_{n,q}(x|\lambda)$ are the $q$-Boole polynomials which are defined by
\begin{equation}
{1\over{1+q(1+t)^{\lambda}}}(1+t)^{x}=\sum_{n=0}^{\infty}Bl_{n,q}(x|\lambda){t^{n}\over{n!}}.
\end{equation}
From (2.1), we can derive the following equation :
\begin{equation}
\int_{\mathbb{Z}_p}\binom{x+\lambda y}{n} d\mu_{-q}(y)={[2]_{q}\over{n!}}Bl_{n,q}(x|\lambda).
\end{equation}
When $x=0, Bl_{n,q}(\lambda)=Bl_{n,q}(0|\lambda)$ are called the $q$-Boole numbers.\\
Now, we observe that
\begin{equation}\begin{split}
(1+t)^{x+\lambda y} &=e^{(x+\lambda y)\log{(1+t)}}\\
&=\sum_{m=0}^{\infty} {(x+\lambda y)^{m}\over{m!}} \bigl(\log(1+t)\bigr)^{m}\\
&=\sum_{m=0}^{\infty}{(x+\lambda y)^{m}\over{m!}}m!\sum_{n=m}^{\infty}S_{1}(n,m){t^{n}\over{n!}}\\
&=\sum_{n=0}^{\infty}\left\{\sum_{m=0}^{n}(x+\lambda y)^{m}S_{1}(n,m)\right\}{t^{n}\over{n!}}.
\end{split}\end{equation}
The $q$-Euler polynomials are defined by the generating function to be
\begin{equation}
{[2]_{q}\over{qe^{t}+1}}e^{xt}=\sum_{n=0}^{\infty}E_{n,q}(x){t^{n}\over{n!}}.
\end{equation}
Note that $\lim_{q\rightarrow 1}{E_{n,q}(x)}=E_{n}(x)$.\\
When $x=0, E_{n,q}=E_{n,q}(0)$ are called the $q$-Euler numbers.
By (1.3), we easily get
\begin{equation}\begin{split}
\int_{\mathbb{Z}_p}e^{(x+y)t}d\mu_{-q}(y) &={[2]_{q}\over{qe^{t}+1}}e^{xt}\\
&=\sum_{n=0}^{\infty}E_{n,q}(x){t^{n}\over{n!}}.
\end{split}\end{equation}
Thus, by (2.6), we get
\begin{equation}
\int_{\mathbb{Z}_p}(x+y)^{n}d\mu_{-q}(y)=E_{n,q}(x),(n\geq0).
\end{equation}
From (2.1), (2.4) and (2.7), we have
\begin{equation}\begin{split}
\int_{\mathbb{Z}_p}(1+t)^{x+\lambda y}d\mu_{-q}(y) &=\sum_{n=0}^{\infty}\left\{\sum_{m=0}^{n}\int_{\mathbb{Z}_p}(x+\lambda y)^{m}d\mu_{-q}(y)S_{1}(n,m)\right\}{t^{n}\over{n!}}\\
&=\sum_{n=0}^{\infty}\left\{\sum_{m=0}^{n}\lambda^{m}E_{m,q}\left({x\over{\lambda}}\right)S_{1}(n,m)\right\}{t^{n}\over{n!}}.
\end{split}\end{equation}
Therefore, by (2.1), (2.3) and (2.8), we obtain the following theorem.

\begin{thm}\label{THEOREM 2.1} For $n\geq0$, we have
\begin{equation*}
Bl_{n,q}(x|\lambda)={1\over{[2]_{q}}}\sum_{m=0}^{n}\lambda^{m}E_{m,q}\left({x\over{\lambda}}\right)S_{1}(n,m),
\end{equation*}
and
\begin{equation*}
\int_{\mathbb{Z}_p}\binom{x+\lambda y}{n}d\mu_{-q}(y)={[2]_{q}\over{n!}}Bl_{n,q}(x|\lambda).
\end{equation*}
\end{thm}

From (2.3), we note that
\begin{equation*}
Bl_{n,q}(x|\lambda)={1\over{[2]_{q}}}\int_{\mathbb{Z}_p}(x+\lambda y)_{n}d\mu_{-q}(y).
\end{equation*}
When $\lambda=1$, we have
\begin{equation}
Bl_{n,q}(x|1)={1\over{[2]_{q}}}\int_{\mathbb{Z}_p}(x+y)_{n}d\mu_{-q}(y).
\end{equation}
As is known, $q$-Changhee polynomials are defined by the generating function to be
\begin{equation}
{[2]_{q}\over{[2]_{q}+qt}}(1+t)^{x}=\sum_{n=0}^{\infty}Ch_{n,q}(x){t^{n}\over{n!}}.
\end{equation}
Thus, by (2.10), we get
\begin{equation}
\int_{\mathbb{Z}_p}(1+t)^{x+y}d\mu_{-q}(y)={[2]_{q}\over{[2]_{q}+qt}}(1+t)^{x}=\sum_{n=0}^{\infty}Ch_{n,q}(x){t^{n}\over{n!}}.
\end{equation}
From (2.11), we have
\begin{equation}\begin{split}
&\int_{\mathbb{Z}_p}(x+y)_{n}d\mu_{-q}(y)=Ch_{n,q}(x),\\
&\textnormal{where}(x)_{n}=x(x-1)\cdots(x-n+1).
\end{split}\end{equation}
By (2.9) and (2.12), we get
\begin{equation}
Bl_{n,q}(x|1)={1\over{[2]_{q}}}Ch_{n,q}(x).
\end{equation}
By replacing $t$ by $e^{t}-1$ in (2.2), we see that
\begin{equation}\begin{split}
{[2]_{q}\over{qe^{\lambda t}+1}}e^{xt} &=[2]_{q}\sum_{n=0}^{\infty}Bl_{n,q}(x|\lambda){1\over{n!}}(e^{t}-1)^{n}\\
&=[2]_{q}\sum_{n=0}^{\infty}Bl_{n,q}(x|\lambda)\sum_{m=n}^{\infty}S_{2}(m,n){t^{m}\over{m!}}\\
&=\sum_{m=0}^{\infty}\sum_{n=0}^{m}[2]_{q}Bl_{n,q}(x|\lambda)S_{2}(m,n){t^{m}\over{m!}},
\end{split}\end{equation}
and
\begin{equation}\begin{split}
{[2]_{q}\over{qe^{\lambda t}+1}}e^{xt} &={[2]_{q}\over{qe^{\lambda t}+1}}e^{\left({x\over{\lambda}}\right)\lambda t}\\
&=\sum_{m=0}^{\infty}E_{m,q}\left({x\over{\lambda}}\right)\lambda^{m}{t^{m}\over{m!}}.
\end{split}\end{equation}
Therefore, by (2.14) and (2.15), we obtain the following theorem.

\begin{thm}\label{THEOREM 2.2} For $m\geq0$, we have
\begin{equation*}
\sum_{n=0}^{m}Bl_{n,q}(x|\lambda)S_{2}(m,n)={1\over{[2]_{q}}}E_{m,q}\left({x\over{\lambda}}\right)\lambda^{m}.
\end{equation*}
\end{thm}

Let us define the $q$-Boole numbers of the first kind with order $k(\in\mathbb{N})$ as follows :
\begin{equation}
[2]_{q}^{k}Bl_{n,q}^{(k)}(\lambda)=\int_{\mathbb{Z}_p}\cdots\int_{\mathbb{Z}_p}(\lambda(x_{1}+\cdots+x_{k}))_{n}d\mu_{-q}(x_{1})\cdots d\mu_{-q}(x_{k}),(n\geq0).
\end{equation}
Thus, by (2.16), we see that
\begin{equation}\begin{split}
[2]_{q}^{k}\sum_{n=0}^{\infty}Bl_{n,q}^{(k)}(\lambda){t^{n}\over{n!}} &=\int_{\mathbb{Z}_p}\cdots\int_{\mathbb{Z}_p}\sum_{n=0}^{\infty}\binom{\lambda(x_{1}+\cdots+x_{k})}{n}t^{n}d\mu_{-q}(x_{1})\cdots d\mu_{-q}(x_{k})\\
&=\int_{\mathbb{Z}_p}\cdots\int_{\mathbb{Z}_p}(1+t)^{\lambda(x_{1}+\cdots+x_{k})}d\mu_{-q}(x_{1})\cdots d\mu_{-q}(x_{k})\\
&=\left({1+q\over{1+q(1+t)^{\lambda}}}\right)^{k}\\
&=[2]_{q}^{k}\sum_{n=0}^{\infty}\left(\sum_{l_{1}+\cdots+l_{k}=n}\binom{n}{l_{1},\cdots,l_{k}}Bl_{l_{1,q}}\cdots Bl_{l_{k,q}}\right){t^{n}\over{n!}}.
\end{split}\end{equation}
Therefore, by (2.17), we obtain the following corollary.

\begin{cor}\label{COROLLARY 2.3} For $n\geq0$, we have
\begin{equation*}
Bl_{n,q}^{(k)}=\sum_{l_{1}+\cdots+l_{k}=n}\binom{n}{l_{1},\cdots,l_{k}}Bl_{l_{1,q}}\cdots Bl_{l_{k,q}}.
\end{equation*}
\end{cor}

The $q$-Euler polynomials of order $k$ are defined by the generating function to be
\begin{equation}\begin{split}
&\int_{\mathbb{Z}_p}\cdots\int_{\mathbb{Z}_p}e^{(x_{1}+\cdots+x_{k}+x)t}d\mu_{-q}(x_{1})\cdots d\mu_{-q}(x_{k})\\
&=\left({[2]_{q}\over{qe^{t}+1}}\right)^{k}e^{xt}=\sum_{n=0}^{\infty}E_{n,q}^{(k)}(x){t^{n}\over{n!}}.
\end{split}\end{equation}
Thus, by (2.18), we get
\begin{equation*}
\int_{\mathbb{Z}_p}\cdots\int_{\mathbb{Z}_p}(x_{1}+\cdots+x_{k}+x)^{n}d\mu_{-q}(x_{1})\cdots d\mu_{-q}(x_{k})=E_{n,q}^{(k)}(x).
\end{equation*}
When $x=0, E_{n,q}^{(k)}=E_{n,q}^{(k)}(0)$ are called the $q$-Euler numbers of order $k$.\\
From (2.16), we note that
\begin{equation}\begin{split}
[2]_{q}^{k}Bl_{n,q}^{(k)}(\lambda) &=\int_{\mathbb{Z}_p}\cdots\int_{\mathbb{Z}_p}(\lambda(x_{1}+\cdots+x_{k}))_{n}d\mu_{-q}(x_{1})\cdots d\mu_{-q}(x_{k})\\
&=\sum_{l=0}^{n}S_{1}(n,l)\int_{\mathbb{Z}_p}\cdots\int_{\mathbb{Z}_p}\lambda^{l}(x_{1}+\cdots+x_{k})^{l}d\mu_{-q}(x_{1})\cdots d\mu_{-q}(x_{k})\\
&=\sum_{l=0}^{n}S_{1}(n,l)\lambda^{l}E_{l,q}^{(k)}.
\end{split}\end{equation}
Therefore, by (2.19), we obtain the following theorem.

\begin{thm}\label{THEOREM 2.4} For $n\geq0$, we have
\begin{equation*}
Bl_{n,q}^{(k)}(\lambda)={1\over{[2]_{q}^{k}}}\sum_{l=0}^{n}S_{1}(n,l)\lambda^{l}E_{l,q}^{(k)}.
\end{equation*}
\end{thm}

By replacing $t$ by $e^{t}-1$ in (2.17), we get
\begin{equation}\begin{split}
[2]_{q}^{k}\sum_{n=0}^{\infty}Bl_{n,q}^{(k)}(\lambda){1\over{n!}}(e^{t}-1)^n&=\left({[2]_{q}\over{qe^{\lambda t}+1}}\right)^{k}\\
&=\sum_{m=0}^{\infty}E_{m,q}^{(k)}\lambda^{m}\frac{t^{m}}{m!},
\end{split}\end{equation}
and
\begin{equation}\begin{split}
[2]_{q}^{k}\sum_{n=0}^{\infty}Bl_{n,q}^{(k)}(\lambda){1\over{n!}}(e^{t}-1)^n&=[2]_{q}^{k}\sum_{n=0}^{\infty}Bl_{n,q}^{(k)}(\lambda)\sum_{m=n}^{\infty}S_{2}(m,n){t^{m}\over{m!}}\\
&=[2]_{q}^{k}\sum_{m=0}^{\infty}\left\{\sum_{n=0}^{m}Bl_{n,q}^{(k)}(\lambda)S_{2}(m,n)\right\}{t^{m}\over{m!}}.
\end{split}\end{equation}
Therefore, by (2.20) and (2.21), we obtain the following theorem.

\begin{thm}\label{THEOREM 2.5} For $m\geq0$, we have
\begin{equation*}
\sum_{n=0}^{m}Bl_{n,q}^{(k)}(\lambda)S_{2}(m,n)={1\over{[2]_{q}^{k}}}E_{m,q}^{(k)}\lambda^{m}.
\end{equation*}
\end{thm}

Let us define the higher-order $q$-Boole polynomials of the first kind as \\
follows :
\begin{equation}\begin{split}
&[2]_{q}^{k}Bl_{n,q}^{(k)}(x|\lambda)=\int_{\mathbb{Z}_p}\cdots\int_{\mathbb{Z}_p}(\lambda x_{1}+\cdots+\lambda x_{k}+x)_{n}d\mu_{-q}(x_{1})\cdots d\mu_{-q}(x_{k}),\\
&\textnormal{where} \ n\geq0 \ \textnormal{and} \ k\in\mathbb{N}.
\end{split}\end{equation}
From (2.22), we can derive the generating function of the higher-order $q$-Boole polynomials of the first kind as follows :
\begin{equation}\begin{split}
[2]_{q}^{k}\sum_{n=0}^{\infty}Bl_{n,q}^{(k)}(x|\lambda){t^{n}\over{n!}}&=\int_{\mathbb{Z}_p}\cdots\int_{\mathbb{Z}_p}(1+t)^{\lambda x_{1}+\cdots+\lambda x_{k}+x}d\mu_{-q}(x_{1})\cdots d\mu_{-q}(x_{k})\\
&=\left({[2]_{q}\over{1+q(1+t)^{\lambda}}}\right)^{k}(1+t)^{x}
\end{split}\end{equation}
By (2.17), we easily get
\begin{equation}\begin{split}
\left({[2]_{q}\over{1+q(1+t)^{\lambda}}}\right)^{k}(1+t)^{x}&=[2]_{q}^{k}\left(\sum_{l=0}^{\infty}Bl_{l,q}^{(k)}(\lambda){t^{l}\over{l!}}\right)\left(\sum_{m=0}^{\infty}m!\binom{x}{m}{t^{m}\over{m!}}\right)\\
&=[2]_{q}^{k}\sum_{n=0}^{\infty}\left(\sum_{m=0}^{n}m!\binom{x}{m}{n!\over{m!(n-m)!}}Bl_{n-m,q}^{(k)}(\lambda)\right){t^{n}\over{n!}}\\
&=[2]_{q}^{k}\sum_{n=0}^{\infty}\left(\sum_{m=0}^{n}m!\binom{x}{m}\binom{n}{m}Bl_{n-m,q}^{(k)}(\lambda)\right){t^{n}\over{n!}}.
\end{split}\end{equation}
Therefore, by (2.23) and (2.24), we obtain the following theorem.

\begin{thm}\label{THEOREM 2.6} For $n\geq0$, we have
\begin{equation*}
Bl_{n,q}^{(k)}(x|\lambda)=\sum_{m=0}^{n}\binom{n}{m}Bl_{n-m,q}^{(k)}(\lambda)(x)_{m}.
\end{equation*}
\end{thm}

Replacing $t$ by $e^{t}-1$ in (2.23), we have
\begin{equation}\begin{split}
[2]_{q}^{k}\sum_{n=0}^{\infty}Bl_{n,q}^{(k)}(x|\lambda){(e^{t}-1)^{n}\over{n!}}&=\left({[2]_{q}\over{1+qe^{\lambda t}}}\right)^{k}e^{xt}\\
&=\sum_{m=0}^{\infty}E_{m,q}^{(k)}\left({x\over{\lambda}}\right)\lambda^{m}{t^{m}\over{m!}},
\end{split}\end{equation}
and
\begin{equation}
[2]_{q}^{k}\sum_{n=0}^{\infty}Bl_{n,q}^{(k)}(x|\lambda){(e^{t}-1)^{n}\over{n!}}=[2]_{q}^{k}\sum_{m=0}^{\infty}\left(\sum_{n=0}^{m}Bl_{n,q}^{(k)}(x|\lambda)S_{2}(m,n)\right){t^{m}\over{m!}}.
\end{equation}
Thus, from (2.25) and (2.26), we have the following theorem.

\begin{thm}\label{THEOREM 2.7} For $m\geq0$ and $k\in\mathbb{N}$, we have
\begin{equation*}
\sum_{n=0}^{m}Bl_{n,q}^{(k)}(x|\lambda)S_{2}(m,n)={1\over{[2]_{q}^{k}}}\lambda^{m}E_{m,q}^{(k)}\left({x\over{\lambda}}\right).
\end{equation*}
\end{thm}

From (2.22), we have
\begin{equation}\begin{split}
[2]_{q}^{k}Bl_{n,q}^{(k)}(x|\lambda)&=\int_{\mathbb{Z}_p}\cdots\int_{\mathbb{Z}_p}(\lambda x_{1}+\cdots+\lambda x_{k}+x)_{n}d\mu_{-q}(x_{1})\cdots d\mu_{-q}(x_{k})\\
&=\sum_{l=0}^{n}S_{1}(n,l)\int_{\mathbb{Z}_p}\cdots\int_{\mathbb{Z}_p}(\lambda x_{1}+\cdots+\lambda x_{k}+x)^{l}d\mu_{-q}(x_{1})\cdots d\mu_{-q}(x_{k})\\
&=\sum_{l=0}^{n}S_{1}(n,l)\lambda^{l}E_{l,q}^{(k)}\left({x\over{\lambda}}\right).
\end{split}\end{equation}
Therefore, by (2.27), we obtain the following theorem.

\begin{thm}\label{THEOREM 2.8} For $n\geq0$, $k\in\mathbb{N}$, we have
\begin{equation*}
Bl_{n,q}^{(k)}(x|\lambda)={1\over{[2]_{q}^{k}}}\sum_{l=0}^{n}S_{1}(n,l)\lambda^{l}E_{l,q}^{(k)}\left({x\over{\lambda}}\right).
\end{equation*}
\end{thm}

Now, we consider the $q$-analogue of Boole polynomials of the second kind as follows :
\begin{equation}
\widehat{Bl}_{n,q}(x|\lambda)={1\over{[2]_{q}}}\int_{\mathbb{Z}_p}(-\lambda y+x)_{n}d\mu_{-q}(y),(n\geq0).
\end{equation}
Thus, by (2.28), we get
\begin{equation}\begin{split}
\widehat{Bl}_{n,q}(x|\lambda)&={1\over{[2]_{q}}}\sum_{l=0}^{n}S_{1}(n,l)(-1)^{l}\lambda^{l}\int_{\mathbb{Z}_p}\left(-{x\over{\lambda}}+y\right)^{l}d\mu_{-q}(y)\\
&={1\over{[2]_{q}}}\sum_{l=0}^{n}S_{1}(n,l)(-1)^{l}\lambda^{l}E_{l,q}\left(-{x\over{\lambda}}\right).
\end{split}\end{equation}
When $x=0, \widehat{Bl}_{n,q}(\lambda)=\widehat{Bl}_{n,q}(0|\lambda)$ are called the $q$-Boole numbers of the second kind.
From (2.28), we can derive the generating function of $\widehat{Bl}_{n,q}(x|\lambda)$ as follows:
\begin{equation}\begin{split}
\sum_{n=0}^{\infty}\widehat{Bl}_{n,q}(x|\lambda){t^{n}\over{n!}}&={1\over{[2]_{q}}}\int_{\mathbb{Z}_p}(1+t)^{-\lambda y+x}d\mu_{-q}(y)\\
&={(1+t)^{\lambda}\over{q+(1+t)^{\lambda}}}(1+t)^{x}.
\end{split}\end{equation}
By replacing $t$ by $e^{t}-1$ in (2.30), we get
\begin{equation}\begin{split}
\sum_{n=0}^{\infty} \widehat{Bl}_{n,q}(x|\lambda){(e^{t}-1)^{n}\over{n!}}&={e^{\lambda t}\over{q+e^{\lambda t}}}e^{xt}\\
&={1\over{qe^{-\lambda t}+1}}e^{xt}\\
&={1\over{[2]_{q}}}\sum_{m=0}^{\infty}(-1)^{m}{\lambda}^{m}E_{m,q}(-\frac{x}{\lambda}){t^{m}\over{m!}},
\end{split}\end{equation}
and
\begin{equation}
\sum_{n=0}^{\infty}\widehat{Bl}_{n,q}(x|\lambda){(e^{t}-1)^{n}\over{n!}}=\sum_{m=0}^{\infty}\left(\sum_{n=0}^{m}\widehat{Bl}_{n,q}(x|\lambda)S_{2}(m,n)\right){t^{m}\over{m!}}.
\end{equation}
Therefore, by (2.31) and (2.32), we obtain the following theorem.

\begin{thm}\label{THEOREM 2.9} For $m\geq0$, we have
\begin{equation*}
{(-1)^{m}\lambda^{m}\over{[2]_{q}}}E_{m,q}(-\frac{x}{\lambda})=\sum_{n=0}^{m}\widehat{Bl}_{n,q}(x|\lambda)S_{2}(m,n),
\end{equation*}
and
\begin{equation*}
\widehat{Bl}_{m,q}(x|\lambda)={1\over{[2]_{q}}}\sum_{l=0}^{m}S_{1}(m,l)(-1)^{l}\lambda^{l}E_{l,q}\left(-{x\over{\lambda}}\right).
\end{equation*}
\end{thm}

For $k\in\mathbb{N}$, let us define the $q$-Boole polynomials of the second kind with order $k$ as follows :
\begin{equation}
\widehat{Bl}_{n,q}^{(k)}(x|\lambda)={1\over{[2]_{q}^{k}}}\int_{\mathbb{Z}_p}\cdots\int_{\mathbb{Z}_p}(-(\lambda x_{1}+\cdots+\lambda x_{k})+x)_{n}d\mu_{-q}(x_{1})\cdots d\mu_{-q}(x_{k}).
\end{equation}
Then we have
\begin{equation*}
[2]_{q}^{k}\widehat{Bl}_{n,q}^{(k)}(x|\lambda)=\sum_{l=0}^{n}S_{1}(n,l)\lambda^{l}(-1)^{l}E_{l,q}\left(-{x\over{\lambda}}\right).
\end{equation*}
From (2.33), we can derive the generating function of $\widehat{Bl}_{n,q}^{(k)}(x|\lambda)$ as follows :
\begin{equation}\begin{split}
\sum_{n=0}^{\infty}\widehat{Bl}_{n,q}^{(k)}(x|\lambda){t^{n}\over{n!}}&={1\over{[2]_{q}^{k}}}\int_{\mathbb{Z}_p}\cdots\int_{\mathbb{Z}_p}(1+t)^{-(\lambda x_{1}+\cdots+\lambda x_{k})+x}d\mu_{-q}(x_{1})\cdots d\mu_{-q}(x_{k})\\
&=\left({(1+t)^{\lambda}\over{q+(1+t)^{\lambda}}}\right)^{k}(1+t)^{x}\\
&=\left({1\over{q(1+t)^{-\lambda}+1}}\right)^{k}(1+t)^{x}\\
&=\sum_{n=0}^{\infty}Bl_{n,q}^{(k)}(x|-\lambda){t^{n}\over{n!}}.
\end{split}\end{equation}
Thus, by (2.34), we get
\begin{equation}
\widehat{Bl}_{n,q}^{(k)}(x|\lambda)={Bl}_{n,q}^{(k)}(x|-\lambda),(n\geq0).
\end{equation}
Indeed,
\begin{equation*}\begin{split}
(-1)^{n}[2]_{q}{Bl_{n,q}(x|\lambda)\over{n!}}&=(-1)^{n}\int_{\mathbb{Z}_p}\binom{x+\lambda y}{n}d\mu_{-q}(y)\\
&=\int_{\mathbb{Z}_p}\binom{-y\lambda-x+n-1}{n}d\mu_{-q}(y)\\
&=\sum_{m=0}^{n}\binom{n-1}{n-m}\int_{\mathbb{Z}_p}\binom{-y\lambda-x}{m}d\mu_{-q}(y)\\
&=\sum_{m=1}^{n}{\binom{n-1}{m-1}\over{m!}}m!\int_{\mathbb{Z}_p}\binom{-y\lambda-x}{m}d\mu_{-q}(y)\\
&=[2]_{q}\sum_{m=1}^{n}\binom{n-1}{m-1}{\widehat{Bl}_{m,q}(-x|\lambda)\over{m!}},
\end{split}\end{equation*}
and
\begin{equation*}\begin{split}
(-1)^{n}[2]_{q}{\widehat{Bl}_{n,q}(x|\lambda)\over{n!}}&=\sum_{m=0}^{n}\binom{n-1}{m-1}\int_{\mathbb{Z}_p}\binom{-x+y\lambda}{m}d\mu_{-q}(y)\\
&=[2]_{q}\sum_{m=1}^{n}\binom{n-1}{m-1}{\widehat{Bl}_{m,q}(-x|\lambda)\over{m!}}.
\end{split}\end{equation*}

\end{document}